# PUMA 560 Trajectory Control using NSGA-II Technique With Real Valued Operators

Habiba BENZATER[1], Samira CHOURAQUI[2]

[1]Department of Computer Science, University of Science and Technology of Oran-Mohammed Boudiaf- (USTO-MB) Oran, Algeria
habibamage@gmail.com

[2]Department of Computer Science, University of Science and Technology of Oran-Mohammed Boudiaf- (USTO-MB) Oran, Algeria
s_chouraqui@yahoo.fr

## ABSTRACT

*In the industry, Multi-objectives problems are a big defy and they are also hard to be conquered by conventional methods. For this reason, heuristic algorithms become an executable choice when facing this kind of problems. The main objective of this work is to investigate the use of the Non-dominated Sorting Genetic Algorithm II (NSGA-II) technique using the real valued recombination and the real valued mutation in the tuning of the computed torque controller gains of a PUMA560 arm manipulator. The NSGA-II algorithm with real valued operators searches for the controller gains so that the six Integral of the Absolute Errors (IAE) in joint space are minimized. The implemented model under MATLAB allows an optimization of the Proportional-Derivative computed torque controller parameters while the cost functions and time are simultaneously minimized.. Moreover, experimental results also show that the real valued recombination and the real valued mutation operators can improve the performance of NSGA-II effectively.*

## KEYWORDS DEFY

*PD Computed Torque Control, Intelligent Control, PUMA560 arm manipulator, Multi-objective Optimization, NSGA-II algorithm, real valued recombination.*

## 1. INTRODUCTION

Multi-Objective Optimization (MOO) has been used for a long time in control engineering. The main target in designing controller is to improve the system's performance by obtaining a good trajectory with minimum error [6] [5].

The computed torque controller proposed first in [10] is a robust nonlinear controller named the computed-torque method in [11, 12]. This technique is known to perform well when the robot arm parameters are known fairly accurately. Fortunately, the dynamic formulation of the PUMA 560 manipulator is well known.

In [19], we have used NSGA-II algorithm introduced in [14] to determine the PD computed torque controller parameters for the PUMA560 system using the SBX crossover and polynomial mutation. The results have shown good performance.

In this paper, we change the recombination operators type of the original NSGA-II algorithm (which are the SBX crossover and the polynomial mutation) to show whether it decrease the quality of the approach or not. In this research, we develop an artificial intelligence (AI) automatic computed torque gains tuning scheme using NSGA-II (Elitist Non-dominated Sorting Genetic Algorithm II) algorithm with real valued recombination and real valued mutation, which can automatically adjust the gains parameters during plant operation in a routine way.

This paper is organized as follows:

In section 2, a general description of the multi-objective evolutionary optimization is presented. Simulation loop of the Puma560 system is presented in section 3. The section 4 presents the NSGA-II algorithm with real valued operators used to tune the PD gains of the computed torque controller. The simulation results of the SIMULINK model and the NSGA-II algorithm with real valued operators are shown in section 5. Finally, we present our conclusions in section 6.

## 2. MULTI-OBJECTIVE EVOLUTIONARY OPTIMIZATION:

In multi-objective optimization (MOO) problems in which the designer seeks to optimize simultaneously several objective functions which are usually in conflict with each other. Such problems haven't usually a single optimal solution. The Pareto dominance relation is the method most commonly adopted to compare solutions in multi-objective optimization which, instead of a single optimal solution, leads to a set of non-dominated solutions called Pareto optimal solutions [9].

Multi-Objective Evolutionary Algorithms (MOEA) are suitable to solve multi-objective optimization problems compared to classical approaches since evolutionary algorithms deal with a set of feasible solutions which allows an efficient way to find an approximation of the whole Pareto optimal solutions in a single simulation run of the algorithm [9].

Nowadays, there are many MOEAs have been suggested, such as Genetic algorithm for multi-objective optimization (MOGA) proposed by Fonseca & Fleming in 1993 [21], the Niched Pareto Genetic Algorithm (NPGA) presented by Horn, Nafpliotis, and Goldberg in 1994 [22], the Non-dominated Sorting Genetic Algorithm (NSGA) introduced by Srinivas & Deb in 1995 [18]. These MOEAS adopt the selection mechanisms based on Pareto ranking and fitness sharing to preserve diversity of the population.

After the algorithms mentioned above, MOEAS based on the elitism strategy were presented, such as the NSGA-II algorithm introduced by Deb, Pratap, Agarwal, & Meyarivan in 2002 [14]. NSGA-II algorithm uses non-dominated sorting method, elitism strategy and a crowded comparison operator for maintaining diversity in the population. The Improved NSGA-II Based on a Novel Ranking Scheme presented by Rio G. L. D'Souza, K. Chandra Sekaran, and A. Kandasamy in [23]. In [24], the authors Long WANG, Tong-guang WANG and Yuan LUO presented a wind turbine blade optimization method which is a novel multi-objective optimization algorithm that employed the controlled elitism strategy and the dynamic crowding distance in the NSGA-II algorithm to improve the lateral diversity and the uniform diversity of non-dominated solutions.

## 3. SIMULATION LOOP OF THE PUMA560 ARM MANIPULATOR:

To simulate the behavior of our system we need several blocks like the trajectory generator block and the controller block as illustrates the Figure.1 given in [8] where $q^d, \tau, \ddot{q}, \dot{q}$ and $q$ are respectively the desired joint positions vector calculated by the trajectory generator, the generalized joint force vector, accelerations vector, velocities vector and position vector.

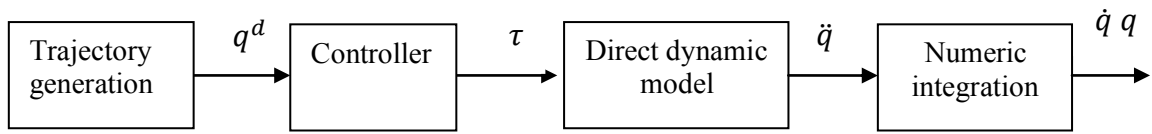

Figure.1 The simulation loop

## 3.1. Dynamic Model of the PUMA560 Arm Manipulator:

The UNIMATION PUMA 560 is a PC controlled, robotic arm that has six revolute joints (or six axes as shown in Figure.2), and each joint is controlled by a DC servo motor and defined by its angle [2] [4].

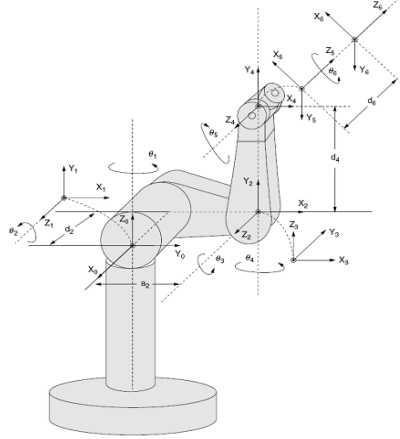

Figure 2.D-H notation for a six-degrees-of-freedom of the PUMA 560 arm manipulator [6]

The dynamic model of the Puma 560 system given in [4] is described in equation (1) as:
$$\tau = A(q)\ddot{q} + B(q)[\dot{q}\dot{q}] + C(q)[\dot{q}^2] + g(q) \qquad (1)$$

Where:

q: is the joint positions (joint angle).
A (q): is the n * n kinetic energy matrix, which is symmetric.
B (q): is the n * n (n-1)/2 matrix of Coriolis torques.
C (q): is the n *n matrix of centrifugal torques.
g (q):is the n-vector of gravity torques;
$\ddot{q}$ : is the n-vector of accelerations.
$\dot{q}$: is the n-vector of joint velocities.
$\tau$ : is the generalized joint force vector.
$[\dot{q}^2]$: $[\dot{q}_1^2, \dot{q}_2^2 \ldots \dot{q}_n^2]$.
$[\dot{q}\dot{q}]$ : $[\dot{q}_1\dot{q}_2, \dot{q}_1\dot{q}_3 \ldots \dot{q}_1\dot{q}_n, \dot{q}_2\dot{q}_3, \dot{q}_2\dot{q}_4, \ldots \dot{q}_{n-2}\dot{q}_n, \dot{q}_{n-1}\dot{q}_n]^T$

The matrices A (q), B (q), C (q) and g (q) can be found in [4].

The direct dynamic model ($\ddot{q}$) given in [8] used to simulate the behavior of the arm manipulator is described in equation (2) as:
$$\ddot{q} = -A^{-1}(q) * [B(q)[\dot{q}\dot{q}] + C(q)[\dot{q}^2] + g(q) - \tau] \qquad (2)$$

## 3.2 Trajectory generator:

Many ways to generate a trajectory both in joint space and Cartesian space as described in [3], [7] and [8]. In our work, we choose the fifth order polynomial given in [8] by the equation (3) to generate motion in the joint space.
$$q(t) = q^i + r(t) * D \qquad 0 \leq t \leq t_f \qquad (3)$$

With:
$$\begin{cases} D = q^f - q^i \\ q^i : \text{Initial joint position} \\ q^f : \text{Final join position} \end{cases} \qquad (4)$$

And:
$$r(t) = 10(t/t_f)^3 - 15(t/t_f)^4 + 6(t/t_f)^5 \qquad (5)$$

### 3.3 The computed torque controller:

Computed torque controller is based on feedback linearization as illustrates figure.3, and computes the required arm torques using the nonlinear feedback control law [5].

We develop this control system in the configuration space, under the assumption that the motion is completely specified, the joint positions and velocities are measurable and that the measurements are not affected with noise.

By using the dynamic equation of the arm given in [8, 3], the computed torque controller is described in equation (6) as:

$$\tau = A(q)\tau' + B(q)[\dot{q}\dot{q}] + C(q)[\dot{q}^2] + g(q) \qquad (6)$$

Where $\qquad \tau' = \ddot{q}_d + k_d * (\dot{q}_d - \dot{q}) + k_p * (q_d - q) \qquad (7)$

$\ddot{q}_d, \dot{q}_d$ and $q_d$ Are respectively the desired joint accelerations, velocities and positions vectors

$\tau$: Generalized joint force vector

$\tau'$: Auxiliary control

$\dot{q}, q$ : Are respectively the joint velocities and positions vectors

$k_p k_d$: Are (n × n) matrices, which are generally diagonal with constant gains on the diagonal

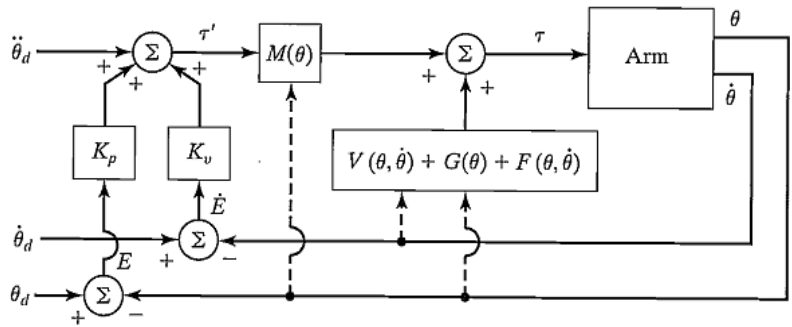

Figure.3 A model-based manipulator-control system [3]

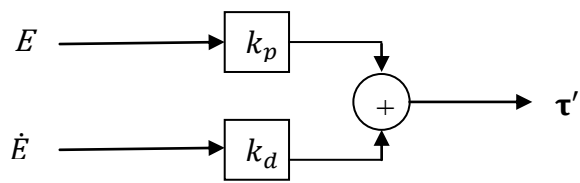

Figure.4 Auxiliary controller internal structure

### 3.5 The SIMULINK diagrams:

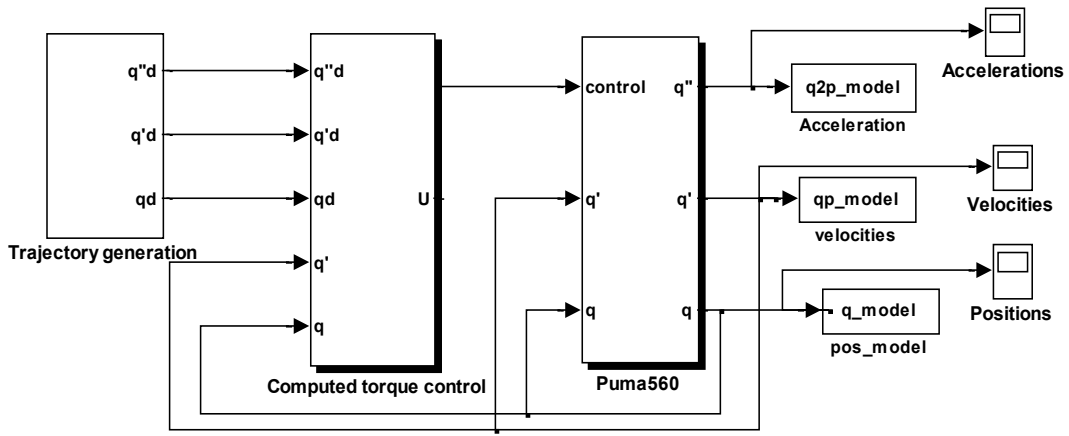

Figure .5 SIMULINK diagram of the simulation loop of the Puma 560 arm manipulator

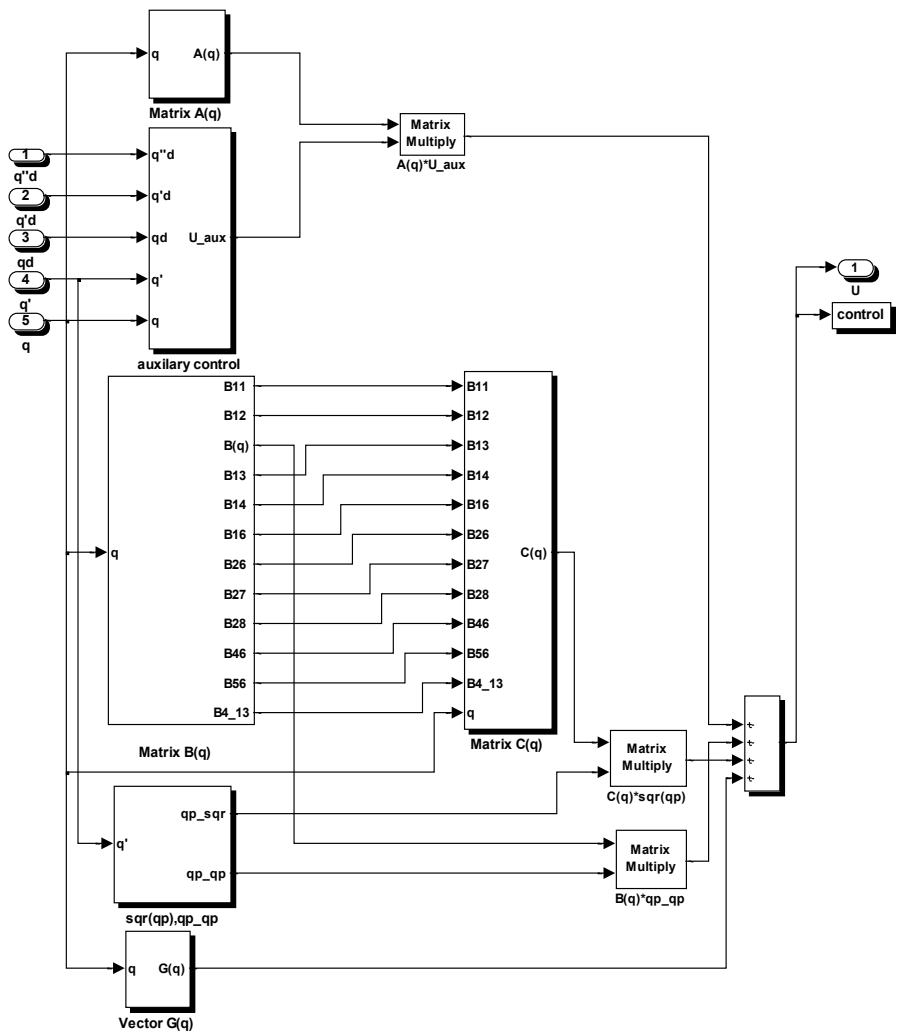

Figure .6 SIMULINK diagram of the computed torque control

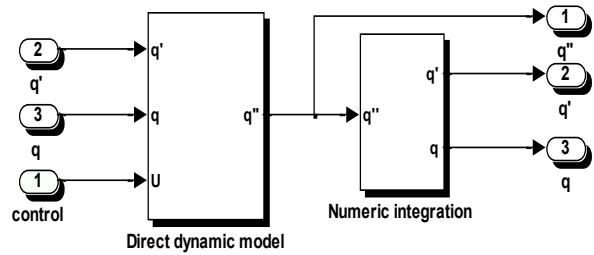

Figure .7 block diagram for q̈, q̇ and q

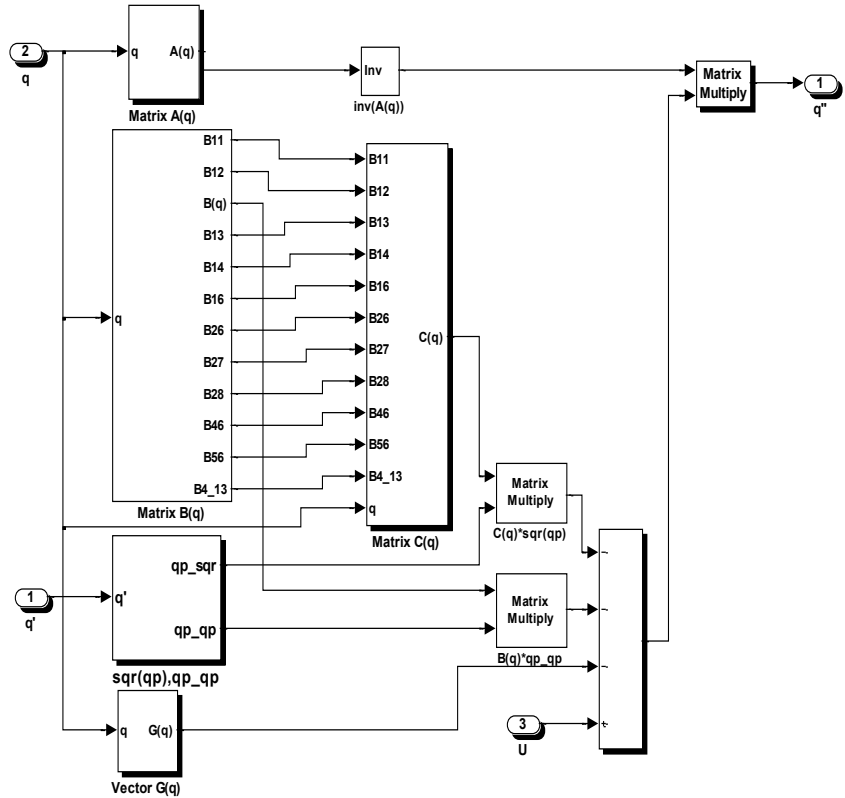

Figure .8 block diagram for q̈

## 4. TUNING PROCEDURE:

### 4.1 NSGA-II algorithm with real valued operators:

Non-dominated Sorting Genetic Algorithm II (NSGA-II) is the heavily revised version of the Non-dominated Sorting Genetic Algorithm (NSGA) which has been generally criticized for its computational complexity, needs elitism operator and requires to choose a priori a value of the sharing function. The NSGA-II algorithm adopts a more efficient ranking procedure. Also, it estimates the crowding distance of the solutions by computing the average distance of two points on either side of this solution along each of the objectives of the problem. In this section a brief description of the NSGA-II technique with real valued operators is given, followed by its application to tune the computed torque controller gains of the Puma560 arm manipulator.

#### 4.1.1 General Description of NSGA-II Algorithm with Real Valued Operators:

4.1.1.1 Initialization of the population by creating N individuals where N>0.
4.1.1.2 Evaluate the m fitness values of each individual in the population.

4.1.1.3 Form and number the fronts based in Pareto dominance relation. The first front is the non-dominant set in the current population, the individuals in the second front are dominated by the individuals in the first front and the fronts go so on.

4.1.1.4 Assign rank value to each individual based on front in which they belong to. Individuals in first front are given a rank value of 1 and individuals in second are assigned rank value as 2 and so on.

4.1.1.5. Estimate the crowding distance for each individual as follow [14]:
For each front $F_i$, n is the number of individuals.
– initialize the distance as zero for all the individuals i.e. $F_i(d_j) = 0$, where j corresponds to the $j^{th}$ individual in front Fi.
– For each objective function m
*Sort the individuals in front $F_i$ i.e. $I = sort(F_i, m)$.
*Assign infinite distance to boundary values in $I$ i.e. $I(d_1) = \infty$ and $I(d_n) = \infty$.
*For k=2 to (n-1)
$I(d_k) = I(d_k) + \frac{I(k+1).m - I(k-1).m}{f_m^{max} - f_m^{min}}$ where $I(k).m$ is the value of the $m^{th}$ objective function of the $k^{th}$ individual in I

4.1.1.6. Selection of the parents by using binary tournament selection for reproduction: During selection, the NSGA-II uses a crowded -comparison operator that supports solutions with minimum rank and in case of solutions that have the same rank, the solutions with the greatest crowding distance is chosen.

4.1.1.7. Generation of the offsprings: The selected parents generate offsprings using Real Valued Recombination and Real Valued Mutation. The real valued recombination technique is applied for the recombination of individuals with real valued variables [20].

4.1.1.8. Apply the step 4.1.1.2 to new offsprings then steps: 4.1.1.3, 4.1.1.4 and 4.1.1.5 to all individuals (selected parents and new offsprings).

4.1.1.9. Based in non-domination relation and crowding distance, sort the whole population formed of the old population and current offspring and select only the N best individuals, where N is the population size.

4.1.1.10. Back to 4.1.1.6 until a convergence criterion is met. More details on the algorithm are found in [14].

## 4.2. NSGAII with Real Valued Operators for PD computed torque controller tuning of puma560:

Using NSGAII algorithm, a solution is represented as a chromosome (i.e.: a string). By taking a 12 variable string as $[k_{p1}\ k_{p2}\ ...\ k_{p6} k_{d1} k_{d2}\ ...\ k_{d6}]$ for NSGAII, an optimal value can be searched. Since Puma560 contains independent controller for each joint, there are 12 values for gains parameters for six joint controllers. The IAE (Integral of Absolute value of Error) performance index given in equation (8) [15] will be used as the fitness vector. In other word, the method of tuning PD parameters using NSGAII with real valued operators is based in minimizing the IAE performance index.

$$IAE = \sum_{k=1}^{k=l} |e(k)| \qquad (8)$$

The equation (9) describes the error vector $e(k)$ as:

$$e(k) = q_d(t) - q(t) \qquad (9)$$

Where $q_d(k)$ is the desired position vector, $q(k)$ is output position vector and $(k)$ is the system error at $l^{th}$ sampling instant.

The implementation of the NSGAII with real valued operators is presented as follow:

1. Produce N initial individual to form initial population. An individual is presented as $[k_{p1}\ k_{p2}\ ...\ k_{p6} k_{d1} k_{d2} ...\ k_{d6}]$ where $k_{pi}\ k_{dj}$ are set in the range of 0 to 100 ($i = 1..6, j = 1..6$).
2. Evaluate fitness $f_j\ j = 1..6$ ($f_j = IAE_j$) for each individual in the population.
3. Sort the current population based on the non-domination sort and crowding distance.
4. Selection of parents by tournament selection based in crowded -comparison operator.
5. Apply the real valued operators (crossover and mutation) find in [20] to generate offsprings.
6. Apply step 2 to evaluate new individuals then step 3 to sort all individuals (parents and offsprings).
7. Selection of the N best individuals based in crowded -comparison operator to form next population.
8. Return to 4 until a predefined number of generations.

NSGA-II algorithm with real valued operators is employed to tune PD computed torque parameters ($k_p\ k_d$) of control system using the model in equation (2). This tuning is based in minimizing simultaneously six position errors.

Parameters used for simulation of NSGA-II with real valued operators are Crossover probability=0.9, Mutation probability=1/12

The generation number and the population size are defined by the user.

## 5. SIMULATION RESULTS:

### 5.1 Result of SIMULINK:

The simulation was implemented in MATLAB/SIMULINK environment. The PD values given in table.1, found empirically in [16] are used to simulate the puma560 model in the SIMULINK environment.

Table.1 PD gains

| Joint | 1 | 2 | 3 | 4 | 5 | 6 |
|---|---|---|---|---|---|---|
| $k_p$ | 700.0 | 1100 | 400.0 | 40.0 | 30.0 | 40.0 |
| $k_d$ | 20.0 | 20.0 | 20.0 | 5.0 | 5.0 | 5.0 |

To generate trajectory in joint space directly, we used as initial and final positions noted respectively $q^i$, $q^f$ found in [17].

$$q^i = (-20°, 60°, -120°, 0°, -30°, 0°)$$
$$q^f = (20°, -60°, -60°, 0°, 30°, 0°)$$

The system is simulated in 1 second and the sampling time is 0.01s

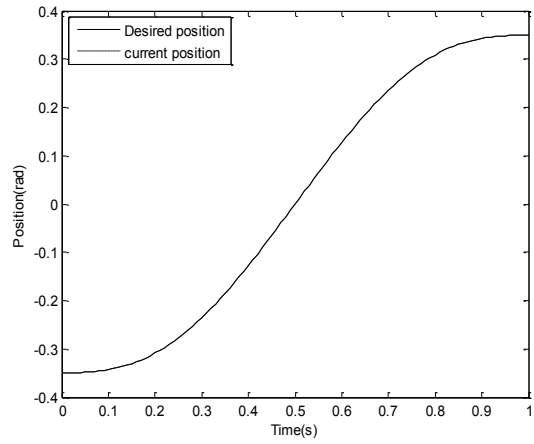

Figure.9 Desired and model joint angles q1

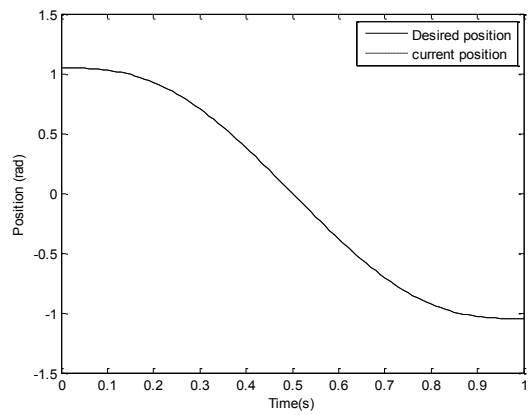

Figure.10 Desired and model joint angles q2

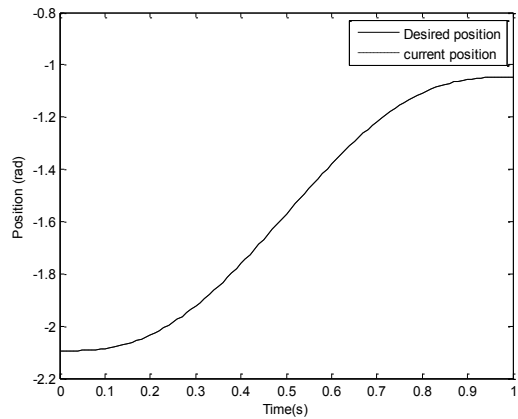

Figure.11 Desired and model joint angles q3

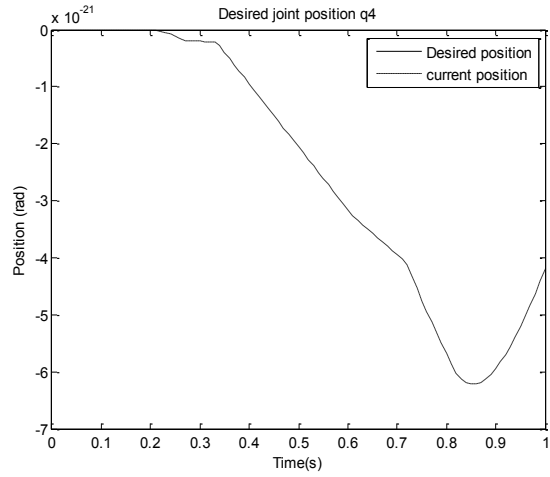

Figure.12 Desired and model joint angles q4

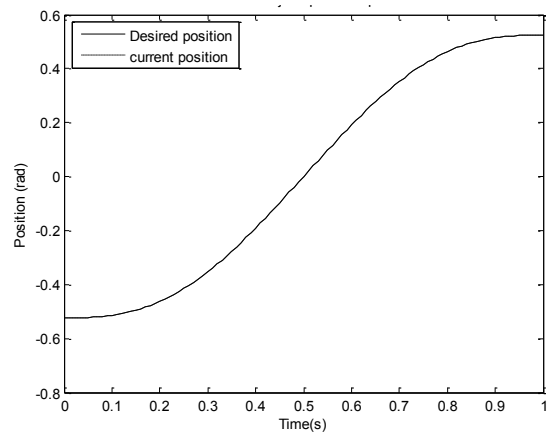

Figure.13 Desired and model joint angles q5

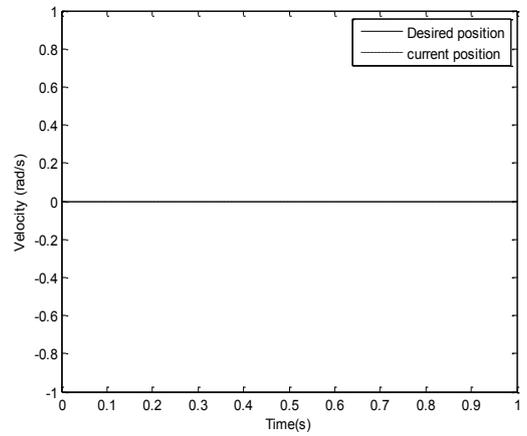

Figure.14 Desired and model joint angles q6

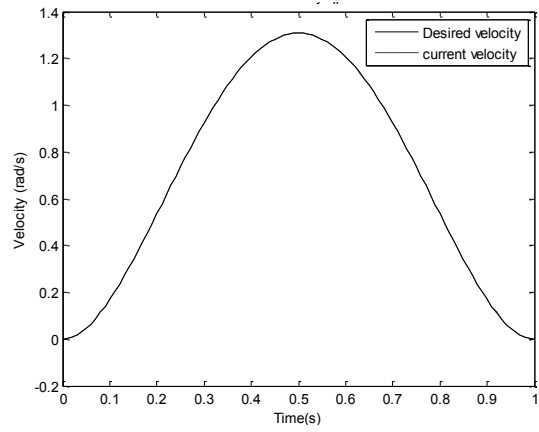

Figure.15 Desired and model joint velocity qp1

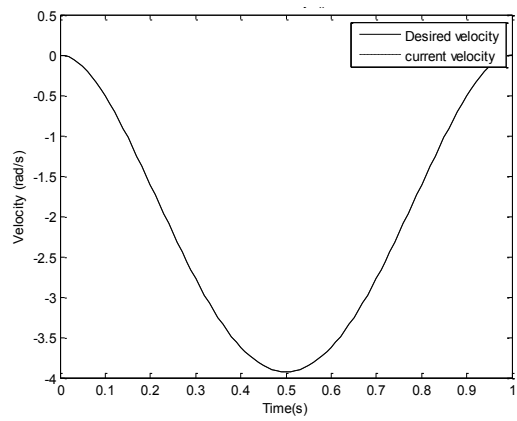

Figure.16 Desired and model joint velocity qp2

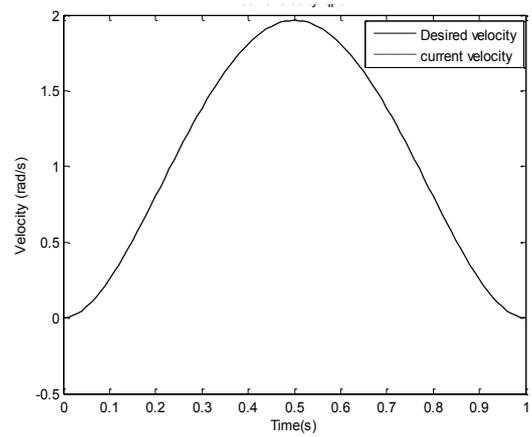

Figure.17 Desired and model joint velocity qp3

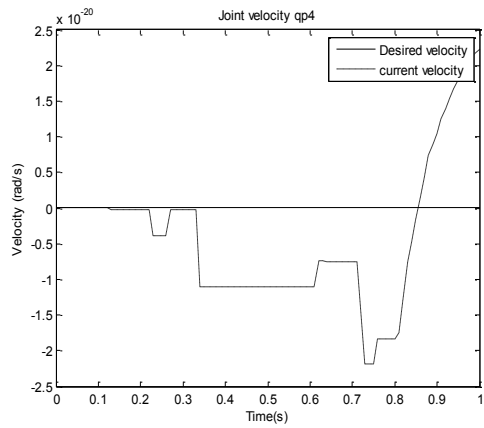

Figure.18 Desired and model joint velocity qp4

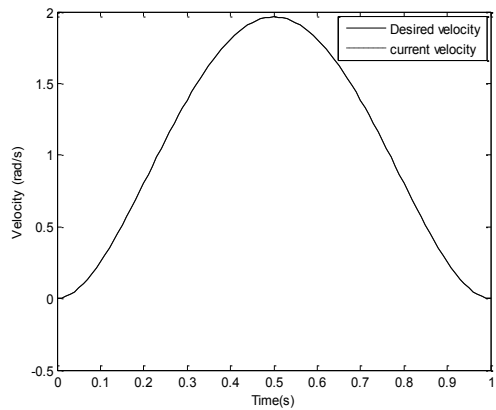

Figure.19 Desired and model joint velocity qp5

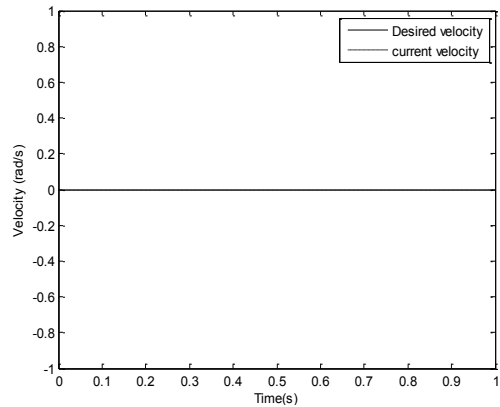

Figure.20 Desired and model joint velocity qp6

## 5.2 Results of NSGA-II Algorithm with real valued operators using for PD Computed Torque Tuning:

In the experiments, the MATLAB R2012a environment has been used for implementing and running the NSGAII algorithm with real valued operators.

As a test, the NSGAII algorithm with real valued operators has been configured as follows: Population size: 2, Generations: 3. Results obtained are: Results are shown in Figure.21 to Figure.32 with an elapsed time 0.3341 s.

- The joint angles (joint positions $q_j \ j = 1,2,3,4,5,6$ ) compared with the desired joint angles values shown in figure.21 to figure.26.
- The joint velocities ($q\dot{p}_j$) compared with the desired joint velocities values shown in figure.27 to figure.32.

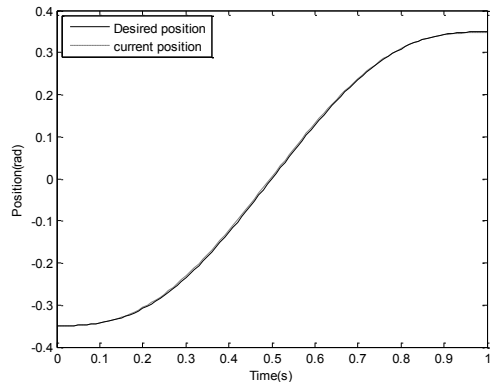

Figure.21 Desired and current joint angles q1

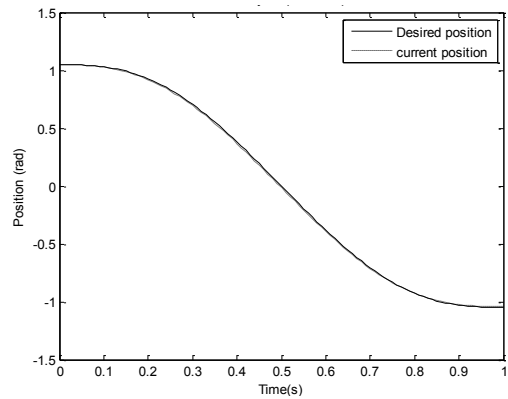

Figure.22 Desired and current joint angles q2

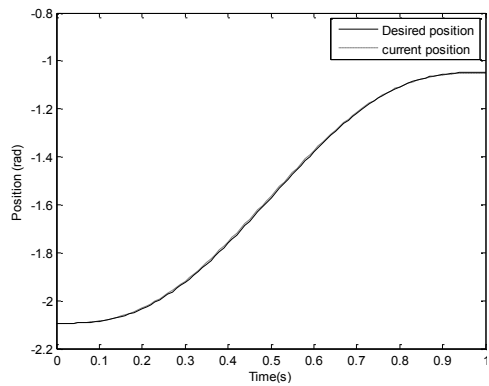

Figure.23 Desired and current joint angles q3

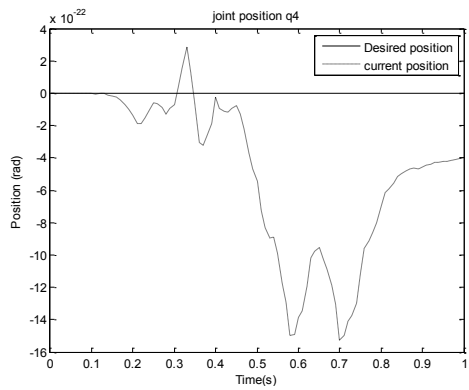
Figure.24 Desired and current joint angles q4

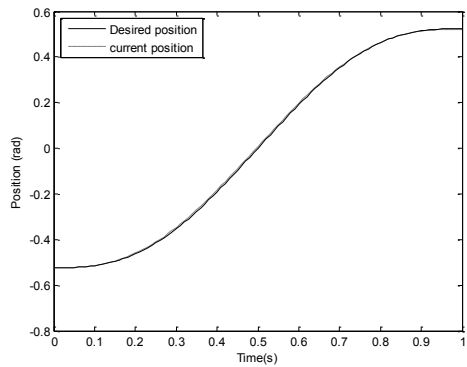
Figure.25 Desired and current joint angles q5

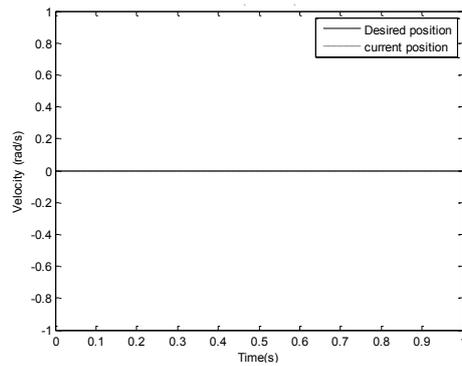
Figure.26 Desired and current joint angles q6

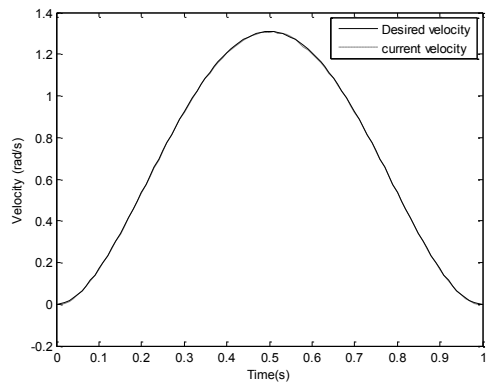
Figure.27 Desired and current joint velocities qp1

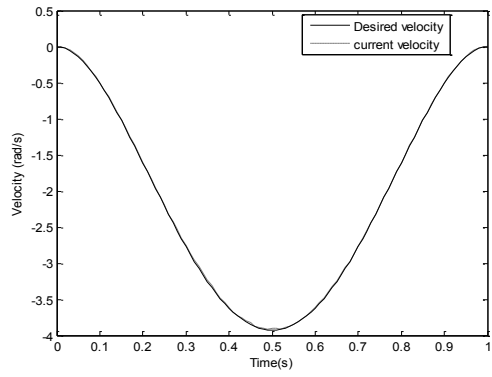

Figure.28 Desired and current joint velocities qp2

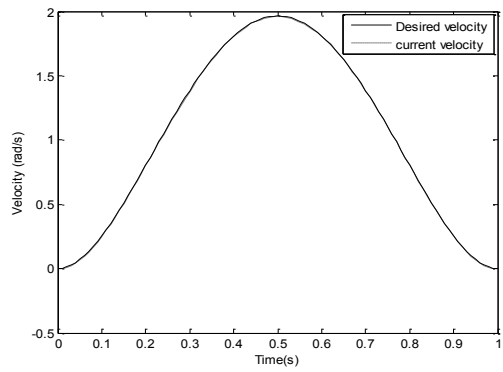

Figure.29 Desired and current joint velocities qp3

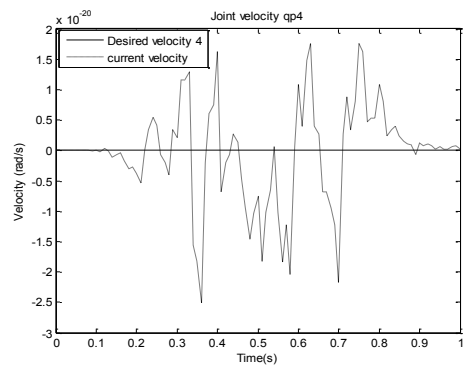

Figure.30 Desired and current joint velocities qp4

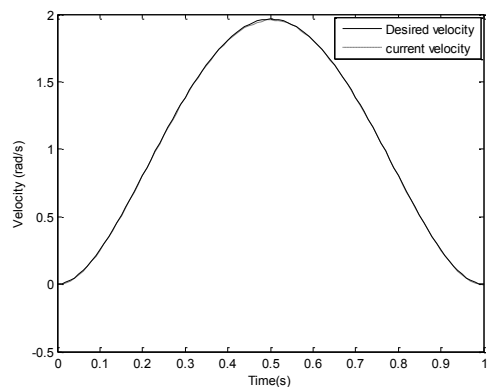

Figure.31 Desired and current joint velocities qp5

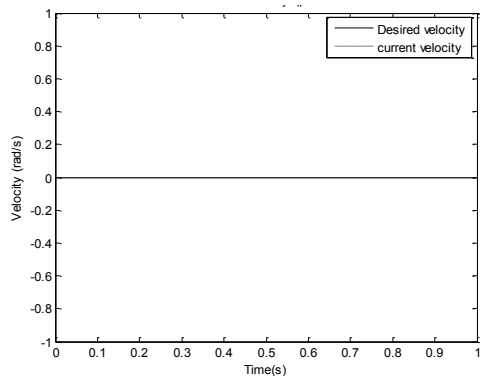

Figure.32 Desired and current joint velocities qp6

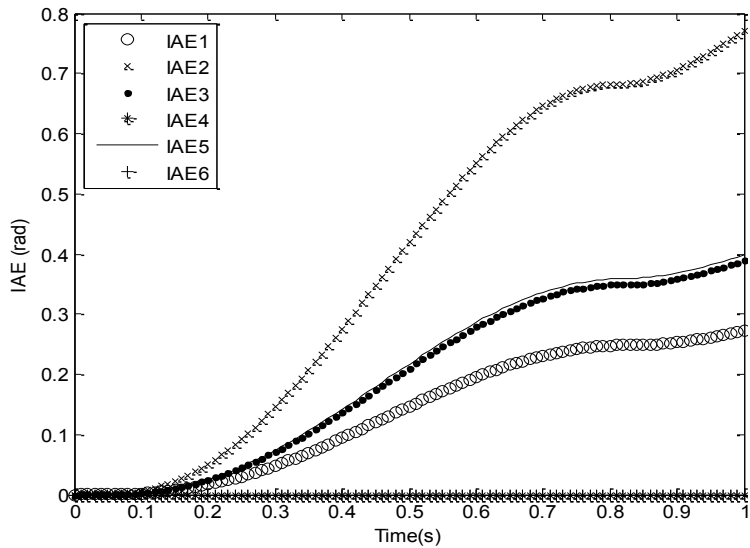

Figure33. The IAE (fitness values) of the joint positions

It can be seen in Figure.21 to Figure.32 that the current joint angles curves resulting from NSGA-II algorithm with real valued converge to the desired joint angles curves. Also, it can be seen that the current joint velocities curves resulting from NSGA-II algorithm with real valued operators converge to the desired joint velocities curves.

It can be observed in figure.9 to figure.20 that the results of SIMULINK model simulations where the gains are tuned empirically are very good.

NSGA-II with real valued operators has proved its effectiveness to regulate the PD-computed torque controller gains of the Puma560 system.

The Figure.34 shows the absolute error of the joint angles obtained by minimizing the IAE errors shown in Figure.33

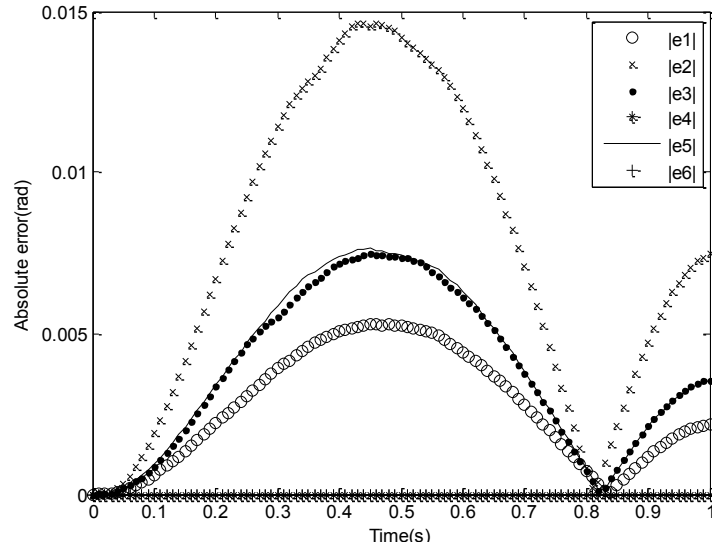

Figure.34 The absolute error of the joint positions

## 6. CONCLUSION

This work is focused on the application of the NSGA-II algorithm using real valued operators to automatically tune the PD-computed torque controller parameters of the PUMA560 arm manipulator. The approach was simulated on Matlab environment to search for the best combination of PD gains so that the Integral of the absolute errors vector in joint space is minimized. It can be concluded from the results that NSGA-II algorithm with real valued operators is showing good performance. As compared to NSGA-II algorithm with SBX crossover and polynomial mutation, it can be concluded that the NSGA-II algorithm with real valued crossover and real valued mutation is one of the recent and efficient optimization tools.

**Authors**

**Habiba Benzater** was born in Algeria and received her engineer degree in computer science from University of Science and Technology of Oran USTO-MB, Algeria in 2009. She is working toward the Magister degree in computer science at University of Science and Technology of Oran USTO-MB, Algeria. . Her current area of research includes optimization and Artificial intelligence.

**SAMIRA CHOURAQUI** received the Ph.D. degree in Computer Science from the University of Science and Technology of Oran-Mohammed Boudiaf Oran, Algeria in 2009. She is a member of the Laboratory of modeling and simulation of the industrial systems. Her current research interests include Artificial Intelligence, Optimization, Fuzzy Control and Neural Networks.

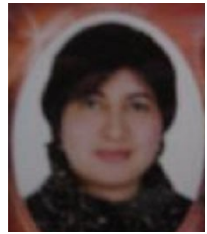